\documentclass[twoside, 11pt]{article}
\usepackage{amssymb, amsmath, mathrsfs, amsthm}
\usepackage{graphicx}
\usepackage{color}
\usepackage{float, caption, subcaption}

\input{epsf}

\newcommand{\bd}{\begin{description}}
\newcommand{\ed}{\end{description}}
\newcommand{\bi}{\begin{itemize}}
\newcommand{\ei}{\end{itemize}}
\newcommand{\be}{\begin{enumerate}}
\newcommand{\ee}{\end{enumerate}}
\newcommand{\beq}{\begin{equation}}
\newcommand{\eeq}{\end{equation}}
\newcommand{\beqs}{\begin{eqnarray*}}
\newcommand{\eeqs}{\end{eqnarray*}}

\definecolor{DarkGreen}{rgb}{0.2, 0.6, 0.3}

\newcommand{\labelz}[1]{\label{#1}}

\newtheorem{theorem}{Theorem}

\newtheorem{lemma}{Lemma}

\newtheorem{case}{Case}
\newtheorem{subcase}{Subcase}[case]
\newtheorem{claim}{Claim}
\newtheorem{fact}{Fact}

\begin{document}
\title{\textbf{Ramsey and Gallai-Ramsey number for wheels\footnote{Supported by the National Science Foundation of China
        (Nos. 11601254, 11551001, 11161037, and 11461054) and the Science
        Found of Qinghai Province (Nos. 2016-ZJ-948Q, and 2014-ZJ-907).}}}

\author{
Yaping Mao\footnote{School of Mathematics and Statistics, Qinghai
Normal University, Xining, Qinghai 810008, China. {\tt
maoyaping@ymail.com}} \footnote{Academy of Plateau Science and Sustainability, Xining, Qinghai 810008, China}, \ \ Zhao Wang\footnote{College of Science,
China Jiliang University, Hangzhou 310018, China. {\tt
wangzhao@mail.bnu.edu.cn}}, \\  Colton Magnant\footnote{Department
of Mathematics, Clayton State University, Morrow, GA, 30260, USA.
{\tt dr.colton.magnant@gmail.com}} \footnotemark[3], \ \  Ingo
Schiermeyer\footnote{Technische Universit{\"a}t Bergakademie
Freiberg, Institut f{\"u}r Diskrete Mathematik und Algebra, 09596
Freiberg, Germany. {\tt Ingo.Schiermeyer@tu-freiberg.de}}. }

\date{}
\maketitle

\begin{abstract}
Given a graph $G$ and a positive integer $k$, define the \emph{Gallai-Ramsey number} to be the minimum number of vertices $n$ such that any $k$-edge coloring of $K_n$ contains either a rainbow (all different colored) triangle or a monochromatic copy of $G$. Much like graph Ramsey numbers, Gallai-Ramsey numbers have gained a reputation as being very difficult to compute in general. As yet, still only precious few sharp results are known. In this paper, we obtain bounds on the Gallai-Ramsey number for wheels and the exact value for the wheel on $5$ vertices.
\end{abstract}

\section{Introduction}
In this work, we consider only edge-colorings of graphs. A coloring of a graph is called \emph{rainbow} if no two edges have the same color.

Colorings of complete graphs that contain no rainbow triangle have interesting and somewhat surprising structure. In 1967, Gallai \cite{MR0221974} first examined this structure under the guise of transitive orientations of graphs. His seminal result in the area was reproven in \cite{MR2063371} in the terminology of graphs and can also be traced to \cite{MR1464337}. For the following statement, a trivial partition is a partition into only one part.

\begin{theorem}[\cite{MR1464337,MR0221974,MR2063371}]\label{Thm:G-Part}
In any coloring of a complete graph containing no rainbow triangle, there exists a nontrivial partition of the vertices (that is, with at least two parts) such that there are at most two colors on the edges between the parts and only one color on the edges between each pair of parts.
\end{theorem}

We refer to a colored complete graph with no rainbow triangle as a \emph{Gallai-coloring} and the partition provided by Theorem~\ref{Thm:G-Part} as a \emph{Gallai-partition}. The induced subgraph of a Gallai colored complete graph constructed by selecting a single vertex from each part of a Gallai partition is called the \emph{reduced graph} of that partition. By Theorem~\ref{Thm:G-Part}, the reduced graph is a $2$-colored complete graph.

Given two graphs $G$ and $H$, let $R(G, H)$ denote the $2$-color Ramsey number for finding a monochromatic $G$ or $H$, that is, the minimum number of vertices $n$ needed so that every red-blue coloring of $K_{n}$ contains either a red copy of $G$ or a blue copy of $H$. Similarly let $R_{k}(H)$ denote the $k$-color Ramsey number for finding a monochromatic copy of $H$ (in any color), that is the minimum number of vertices $n$ needed so that every $k$-coloring of $K_{n}$ contains a monochromatic copy of $H$. Although the reduced graph of a Gallai partition uses only two colors, the original Gallai-colored complete graph could certainly use more colors. With this in mind, we consider the following generalization of the Ramsey numbers. Given two graphs $G$ and $H$, the \emph{general $k$-colored Gallai-Ramsey number} $gr_k(G:H)$ is defined to be the minimum integer $m$ such that every $k$-coloring of the complete graph on $m$ vertices contains either a rainbow copy of $G$ or a monochromatic copy of $H$. With the additional restriction of forbidding the rainbow copy of $G$, it is clear that $gr_k(G:H)\leq R_k(H)$ for any graph $G$.

Recently, there has been a flurry of activity in the area with an influx of new results and approaches. In particular, the following results were recently obtained for fans.

\begin{theorem}[\cite{Fans}]
$$
gr_k(K_3;F_2)=\begin{cases}
9, &\mbox {\rm if}~k=2;\\[0.2cm]
\frac{83}{2}\cdot 5^{\frac{k-4}{2}}+\frac{1}{2}, &\mbox {\rm if}~k~is~even,~k\geq 4;\\[0.2cm]
4\cdot 5^{\frac{k-1}{2}}+1, &\mbox {\rm if}~k~is~odd.
\end{cases}
$$
\end{theorem}

\begin{theorem}[\cite{Fans}]
For $k\geq 2$,
$$
\begin{cases}
4n\cdot 5^{\frac{k-2}{2}}+1\leq gr_k(K_3;F_n)\leq 10n\cdot 5^{\frac{k-2}{2}}-\frac{5}{2}n+1, &\mbox {\rm if}~k~is~even;\\[0.2cm]
2n\cdot 5^{\frac{k-1}{2}}+1\leq gr_k(K_3;F_n)\leq \frac{9}{2}n\cdot 5^{\frac{k-1}{2}}-\frac{5}{2}n+1, &\mbox {\rm if}~k~is~odd.
\end{cases}
$$
\end{theorem}

Odd cycles were also recently settled completely.

\begin{theorem}[\cite{OddCycles}] 
For integers $\ell \geq 3$ and $k \geq 1$, we have
$$
gr_{k} (K_{3} : C_{2\ell + 1}) = \ell \cdot 2^{k} + 1.
$$
\end{theorem}

In this work, we consider the Gallai-Ramsey numbers for finding either a rainbow triangle or monochromatic wheel. Let $W_n$ be a wheel of order $n$, that is, $W_{n}=K_{1}\vee C_{n-1}$ where $C_{n - 1}$ is the cycle on $n - 1$ vertices.
\begin{theorem}{\upshape \cite{MR1670625}}\label{Thm:Ramsey}
$$
\begin{array}{l}
\displaystyle (1) ~ R(W_{5}, W_{5})=15;\\
\displaystyle (2) ~ R(W_{6}, W_{6})=17.
\end{array}
$$
\end{theorem}

As far as we are aware, for $n\geq 7$, the classical diagonal Ramsey number for the wheel is yet unknown. We give upper and lower bounds for classical Ramsey number of the general wheel $W_{n}$ in Section~\ref{Sec:Ram}.
\begin{theorem}\label{Thm:RamseyWn}
For $k \geq 1$ and $n\geq 7$,
$$
\begin{cases}
3n-3\leq R(W_{n},W_{n})\leq 8n-10, & \text{ if $n$ is even;}\\
2n-2\leq R(W_{n},W_{n})\leq 6n-8 & \text{ if $n$ is odd.}
\end{cases}
$$
\end{theorem}

In Section~\ref{Sec:W5}, we obtain the exact value of the Gallai Ramsey number
for $W_5$.
\begin{theorem}\label{Thm:GallaiRamseyW5}
For $k \geq 1$,
$$
gr_{k}(K_{3}:W_{5})=\begin{cases}
5 & \text{ if $k = 1$,}\\
14\cdot5^{\frac{k-2}{2}}+1 & \text{ if $k$ is even,}\\
28\cdot5^{\frac{k-3}{2}}+1 & \text{ if $k \geq 3$ is odd.}
\end{cases}
$$
\end{theorem}


Finally in Section~\ref{Sec:Gen}, we provide general upper and lower bounds on the Gallai-Ramsey numbers for all wheels.

We refer the interested reader to \cite{MR1670625} for a dynamic
survey of small Ramsey numbers and \cite{FMO14} for a dynamic survey
of rainbow generalizations of Ramsey theory, including topics like
Gallai-Ramsey numbers.

\section{Bounds on the Ramsey numbers}\label{Sec:Ram}

First some additional definitions. A cycle $C_k$ of length $k$ is also called a \emph{$k$-cycle}. A path of a graph $G$ is a Hamiltonian path if it contains all the vertices of $G$. A graph $G$ is said to be \emph{pancyclic} if it has $k$-cycles for every $k$ between $3$ and $n$. A vertex of a graph $G$ is \emph{$r$-pancyclic} if it is contained in a $k$-cycle for every $k$ between $r$ and $n$, and $G$ is \emph{vertex $r$-pancyclic} if every vertex is $r$-pancyclic.

Hendry \cite{Hendry} derived the following result.
\begin{lemma}[\cite{Hendry}]\label{lem1}
Let $G$ be a graph of order $n\geq 3$ with $\delta(G)\geq (n +
1)/2$. Then $G$ is vertex pancyclic.
\end{lemma}

\begin{lemma}[\cite{FS1, KaRos, Ros1}]\label{lem2}
For $k \geq 1$,
$$
R(C_{m},C_{n})=\begin{cases}
2n-1, & \text{if $3\leq m\leq n$, $m$ odd,} \\
 & ~ (m,n)\neq (3,3);\\
n-1+m/2, & \text{if $4\leq m\leq n$, $m$ and $n$ even,}\\
~ & (m,n)\neq (3,3);\\
\max\{n-1+m/2,2m-1\}, & \text{if $4\leq m\leq n$,}\\
~ & \text{$m$ even and $n$ odd.}
\end{cases}
$$
\end{lemma}

By the above results, we derive the upper and lower bounds for the Ramsey number of general wheels.
\begin{lemma}\label{Lemma:RamseyWnBounds}
For $k \geq 1$ and $t\geq 3$,
$$
6t+4 \leq R(W_{2t+2},W_{2t+2}) \leq 16t+6,
$$
and
$$
4t + 1 \leq R(W_{2t+1},W_{2t+1}) \leq 12t-2.
$$
\end{lemma}

\begin{proof}
First the even case. For the lower bound, let $G$ be a $2$-edge colored graph obtained from three blue copies of $K_{2t+1}$ by adding all red edges in between them. Clearly, there is neither a red copy of $W_{2t+2}$ nor a blue copy of $W_{2t+2}$ in $G$. Since $|G| = 6t + 3$, this means that $R(W_{2t+2}, W_{2t+2})\geq 6t+4$.

Let $G$ be a $2$-edge colored copy of $K_{16t+6}$ with colors red and blue. For each $v\in V(G)$, let $A_v$ and $B_v$ be the set of vertices incident to $v$ be red and blue edges, respectively. Then for every vertex $v\in V(G)$ such that $|A_v|\geq 8t+3$ or $|B_v|\geq 8t+3$. Without loss of generality, we suppose $|A_v|\geq 8t+3$. For each vertex $u \in A_{v}$, let $D_{u}$ be the set of vertices in $A_{v}$ with blue edges to $u$.  If there is a vertex $u \in A_{v}$ with $|D_u|\geq 4t+1$, then since $R(C_{2t+1},C_{2t+1})=4t+1$ (by Lemma~\ref{lem2}), there exists either a red cycle $C_{2t+1}$ or a blue cycle $C_{2t+1}$ within $D_u$. If it is a red cycle $C_{2t+1}$, then the subgraph induced by $V(C_{2t+1})\cup \{v\}$ is a red copy of $W_{2t+2}$. If it is a blue cycle $C_{2t+1}$, then the subgraph induced by $V(C_{2t+1})\cup \{u\}$ is a blue copy of $W_{2t+2}$. Thus, we may assume that for any $u\in A_v$, $|D_u|\leq 4t$. Then the number of incident red edges to $u$ in $A_v$ is at least $|A_v|-4t-1 \geq \frac{|A_v|+1}{2}$. From Lemma~\ref{lem1}, there is red cycle $C_{2t+1}$ in $A_v$. This cycle together with $v$ is a red copy of $W_{2t+2}$. Thus, we have $R(W_{2t+2},W_{2t+2})\leq 16t+6$.

Now the odd case. Let $G$ be a $2$-edge colored graph obtained from two blue copies of $K_{2t}$ by adding all red edges in between them. Clearly, there is neither a red copy of $W_{2t+1}$ nor a blue copy of $W_{2t+1}$ in $G$. Since $|G| = 2(2t) = 4t$, we have $R(W_{2t+2}, W_{2t+2})\geq 4t + 1$.

Much like the proof of the even case, let $G$ be a $2$-edge colored copy of $K_{16t+6}$ with colors red and blue. For each $v\in V(G)$, let $A_v$ and $B_v$ be the set of vertices incident to $v$ be red and blue edges, respectively. Then for every vertex $v\in V(G)$ such that $|A_v|\geq 6t-1$ or $|B_v|\geq 6t-1$. Without loss of generality, we suppose $|A_v|\geq 6t-1$. For each vertex $u \in A_{v}$, let $D_{u}$ be the set of vertices in $A_{v}$ with blue edges to $u$.  If there is a vertex $u \in A_{v}$ with $|D_u|\geq 3t-1$, then since $R(C_{2t},C_{2t})=3t-1$ (by Lemma~\ref{lem2}), there exists either a red cycle $C_{2t}$ or a blue cycle $C_{2t}$ within $D_u$. If it is a red cycle $C_{2t}$, then the subgraph induced by $V(C_{2t})\cup \{v\}$ is a red copy of $W_{2t+1}$. If it is a blue cycle $C_{2t}$, then the subgraph induced by $V(C_{2t})\cup \{u\}$ is a blue copy of $W_{2t+1}$. Thus, we may assume that for any $u\in A_v$, $|D_u|\leq 3t-2$. Then the number of incident red edges to $u$ in $A_v$ is at least $|A_v|-3t-3 \geq \frac{|A_v|+1}{2}$. From Lemma~\ref{lem1}, there is red cycle $C_{2t}$ in $A_v$. This cycle together with $v$ is a red copy of $W_{2t+1}$. Thus, we have $R(W_{2t+1},W_{2t+1})\leq 12t-2$.
\end{proof}

\section{The Gallai-Ramsey number for $W_{5}$} \label{Sec:W5}

In this section, we give the results for the Gallai Ramsey number of $W_5$ and general wheel $W_n$ for $n \geq 6$.

We first give the lower bound on the Gallai-Ramsey number for $W_{5}$.

\begin{lemma}\label{Lemma:W5Lower}
For $k \geq 2$,
$$
gr_{k}(K_{3} : W_{5}) \geq \begin{cases}
14\cdot 5^{(k - 2)/2} + 1 & \text{if $k$ is even,}\\
28\cdot 5^{(k - 3)/2} + 1 & \text{if $k$ is odd.}
\end{cases}
$$
\end{lemma}

\begin{proof}
We prove this result by inductively constructing a coloring of $K_{n}$ where
$$
n = \begin{cases}
14 \cdot 5^{(k - 2)/2} & \text{if $k$ is even,}\\
28 \cdot 5^{(k - 3)/2} & \text{if $k$ is odd,}
\end{cases}
$$
which contains no rainbow triangle and no monochromatic copy of $W_{5}$. For the base of this induction, let $G_{2}$ be a $2$-colored complete graph on $R(W_{5},W_{5}) - 1 = 14$ vertices containing no monochromatic copy of $W_{5}$. Suppose this coloring uses colors $1$ and $2$.

Suppose we have constructed a coloring of $G_{2i}$ where $i$ is a positive integer and $2i < k$, using the $2i$ colors in the set $[2i]$ and having order $n_{2i} = 14 \cdot 5^{i - 1}$ such that $G_{2i}$ contains no rainbow triangle and no monochromatic copy of $W_{5}$.

If $k = 2i + 1$, we construct $G_{2i + 1}$ by making two copies of $G_{2i}$ and inserting all edges between the copies in color $k$. Then $G_{k}$ certainly contains no rainbow triangle, no monochromatic copy of $W_{5}$, and has order $n = 2 \cdot 14 \cdot 5^{(k - 3)/2} = 28 \cdot 5^{(k - 3)/2}$.

Otherwise suppose $k \geq 2i + 2$. We construct $G_{2i + 2}$ by
making five copies of $G_{2i}$ and inserting edges of colors $2i +
1$ and $2i + 2$ between the copies to form a blow-up of the unique
$2$-colored $K_{5}$ which contains no monochromatic triangle. This
coloring clearly contains no rainbow triangle and, since there is no
monochromatic triangle in either of the two new colors, there can be
no monochromatic copy of $W_{5}$, completing the construction.
\end{proof}

We are now in a position to prove Theorem~\ref{Thm:GallaiRamseyW5}, that is, for $k \geq 1$,
$$
gr_{k}(K_{3}:W_{5})=\begin{cases}
5 & \text{ if $k = 1$,}\\
14\cdot5^{\frac{k-2}{2}}+1 & \text{ if $k$ is even,}\\
28\cdot5^{\frac{k-3}{2}}+1 & \text{ if $k \geq 3$ is odd.}
\end{cases}
$$

\begin{proof}
The lower bound follows from Lemma~\ref{Lemma:W5Lower}. Call a color \emph{wasted} if it induces only a matching and \emph{useful} if there are adjacent edges in the color. Note that in a colored complete graph, in order to avoid a rainbow triangle, all wasted colors must together induce a matching. Let $G$ be a $k$-coloring of a complete graph in which there are only $k'$ colors which induce a subgraph containing adjacent edges. If $k'=0$, then every color is wasted and so every set of three vertices induces a rainbow triangle, clearly a contradiction. For $k' \geq 1$, let $n$ be the order of $G$ where
$$
n = n_{k'} = \begin{cases}
5 & \text{ if $k' = 1$,}\\
14\cdot5^{\frac{k'-2}{2}}+1 & \text{ if $k'$ is even,}\\
28\cdot5^{\frac{k'-3}{2}}+1 & \text{ if $k' \geq 3$ is odd.}
\end{cases}
$$
We now prove the upper bound by induction on $k'$ since $k \leq k'$. If $k'=1$, then $G$ is a coloring of $K_{5}$ in which each color is wasted except one, say color $1$. This means that the subgraph induced by color $1$ is a $K_{5}$ minus a matching, which is a copy of $W_{5}$, a contradiction.

Next suppose $k'=2$, so $n=15$. If $G$ uses exactly $2$ colors, then it follows from the fact that $R(W_5,W_5)=15$ that there is a monochromatic copy of $W_5$ in $G$, a contradiction. Suppose, therefore, that $G$ uses at least $3$ colors. Let red and blue be two useful colors. Since all wasted colors induce a single matching, we may assume all wasted edges are green and let $uv$ be a green edge. To avoid a rainbow triangle, every vertex in $G$ other than $u$ and $v$ have a single color (red or blue) to both $u$ and $v$. This being the case with all green edges, there is a Gallai partition of $G$ with all parts of order at most $2$ consisting of the green edges. Let $A$ be the set of parts with red edges to $\{u,v\}$, and $B$ be the set of parts with blue edges to $\{u,v\}$. In order to avoid a red copy of $W_{5}$, there is no vertex in $A$ with two incident red edges within $A$. This means that the red edges within $A$ form a matching, along with any green edges. If $|A| \geq 5$, then since $A$ contains all blue edges except for possibly a matching of red or green edges, there is a blue copy of $W_{5}$ within $A$. Thus, we may assume $|A| \leq 4$ and similarly $|B| \leq 4$. Hence $|G|=|A|+|B|+2 \leq 10<15$, a contradiction.

Suppose $k'\geq 3$. Inductively we suppose the statement is true for all $k'<k$ and consider $k' \doteq k$.

By Theorem~\ref{Thm:G-Part}, there exists a partition of $V(G)$ into parts such that between each pair of parts there is exactly one color and between the parts in general, there are at most two colors (say color $c_{1}$ and $c_{2}$). Consider such a $G$-partition with the smallest number of parts, say $t_1$. Since $R(W_{5}, W_{5})=15$, it follows that $t_1\leq 14$. Let $H_{1}^{1}, H_{2}^{1}, \cdots, H_{t_1}^{1}$ be parts of the $G$-partition.

Suppose $2\leq t_1 \leq 3$. By the minimality of $t_1$, we may assume $t_1=2$. If $|H_{1}^1|=1$, then $|H_{2}^1|\geq 14$ since $n\geq 15$. Without loss of generality, suppose that all edges between $H_{1}^{1}$ and $H_{2}^{1}$ are color $c_{1}$. If there is no edge with color $c_{1}$ in $H_{2}^{1}$, then
$$
|G|=|H_{1}^{1}|+|H_{2}^{1}|\leq 1+[gr_{k^{'}-1}(K_{3}:W_{5})-1]<n,
$$
a contradiction.

We now suppose there are some edges with color $c_{1}$ in $H_{2}^{1}$. Because there is no rainbow triangle or monochromatic copy of $W_{5}$ in $H_{2}^{1}$, by Theorem~\ref{Thm:G-Part}, there exists a partition of $V(H_{2}^{1})$ into parts such that between each pair of parts there is exactly one color and between the parts in general, there are at most two colors. Consider such a $H_{2}^{1}$-partition with the smallest number of parts, say $t_2$, clearly, $2\leq t_2 \leq 14$. Let $H_{1}^{2}, H_{2}^{2}, \cdots, H_{t_2}^{2}$ be parts of the $H_{2}^{1}$-partition.

Suppose $2\leq t_2\leq 3$. By the minimality of $t_2$, we may assume $t_2=2$. If $H_{1}^{2}=1$, then $|H_{2}^{2}|\geq 13$. We suppose that all edges between $H_{1}^{2}$ and $H_{2}^{2}$ are color $c_{1}$. To avoid a $W_{5}$ with color $c_{1}$, there is no $P_{3}$ with color $c_{1}$ within $H_{2}^{2}$, and hence the subgraph induced by color $c_{1}$ is a matching of $H_{2}^{2}$. Then $|H_{2}^{2}|\leq gr_{k'-1}(K_{3}:W_{5})-1$, and so
$$
|G|=|H_{1}^{2}|+|H_{2}^{2}|\leq 2+[gr_{k'-1}(K_{3}:W_{5})-1]<n,
$$
a contradiction. We now suppose that all edges between $H_{1}^{2}$ and $H_{2}^{2}$ are not color $c_{1}$, say $c_2$.

Continue this above process. Then there exists a sequence of vertices $v_{1}, v_{2}, \cdots, v_{s}$ in $G$ such that
\begin{itemize}
\item $H_{1}^{1}=\{v_{1}\}, H_{1}^{2}=\{v_{2}\}, \cdots, H_{1}^{s}=\{v_{s}\}$;
\item Let $H_{1}^{i}, H_{2}^{i}, \cdots, H_{t_i}^{i}$ be parts of the $H_{2}^{i-1}$-partition for each $i \ (1\leq i \leq s)$;
\item $t_1=t_2=\cdots=t_s=2$;
\item The edges from $H_{1}^{i}=\{v_i\}$ to $H_{2}^{i}$ are colored by $c_{i}$ for each $i$ with $1\leq i \leq s$.
\end{itemize}

\begin{claim}\label{Claim1}
$2\leq s\leq 3k'$.
\end{claim}
\begin{proof}
Assume, to the contrary, that $s\geq 3k'+1$. Then there exist $4$
vertices of $v_1,v_2,\ldots,v_s$, say
$v_{i_1},v_{i_2},v_{i_3},v_{i_4}$ ($i_1\leq i_2\leq i_3\leq i_4$),
such that the edges from $v_{i_p}$ to $v_{i_q}$ ($1\leq p\neq q\leq
4$) receives color $c_{i_1}$, the edges from each $v_{i_p} \ (1\leq
p\leq 4)$ to $H_2^s$ receives color $c_{i_1}$. It is clear that
there is a $W_5$ with color $c_{i_1}$, a contradiction.
\end{proof}

Furthermore, we have the following claim by Claim \ref{Claim1}.
\begin{claim}\label{Claim2}
$2\leq s\leq k'$.
\end{claim}
\begin{proof}
Assume, to the contrary, that $s\geq k'+1$. Then there exist at
least $2$ vertices of $v_1,v_2,\ldots,v_s$, say $v_{i_1},v_{i_2}$,
such that the edges from $v_{i_1}$ to $v_{i_2}$ receives color
$c_{i_1}$, the edges from each $v_{i_1}$ to $H_2^s$ receives color
$c_{i_1}$, and the edges from each $v_{i_2}$ to $H_2^s$ receives
color $c_{i_1}$. Then $H_2^s$ contains at least a matching with
color $c_{i_1}$, and hence $|H_2^s|\leq gr_{k'-1}(K_3;W_5)-1$. From
Claim~\ref{Claim1}, we have
$$
|G|=(s-1)+[gr_{k'-1}(K_{3}:W_{5})-1]\leq
(3k'-1)+[gr_{k'-1}(K_{3}:W_{5})-1]<n,
$$
a contradiction.
\end{proof}

If $H_{1}^{s}=1$, then $|H_{2}^{s}|\geq 15-s$. We suppose that all edges between $H_{1}^{s}$ and $H_{2}^{s}$ are color $c_{j} \ (1\leq j\leq s-1)$. To avoid a $W_{5}$ with color $c_{j}$, there is no $P_{3}$ with color $c_{j}$ within $H_{2}^{s}$, and hence the subgraph induced by color $c_{j}$ is a matching of $H_{2}^{s}$. Then $|H_{2}^{s}|\leq gr_{k'-1}(K_{3}:W_{5})-1$, and so
$$
|G|=(s-1)+|H_{1}^{s}|+|H_{2}^{s}|\leq
s+[gr_{k'-1}(K_{3}:W_{5})-1]<n,
$$
a contradiction. We now suppose that all edges between $H_{1}^{s}$ and $H_{2}^{s}$ are not color $c_{j} \ (1\leq j\leq s-1)$, say $c_s$.

Note that $H_{2}^{s}=G-\{v_{1}, v_{2}, \cdots, v_{s}\}$. Then, by Theorem \ref{Thm:G-Part}, we see that $H_{2}^{s}$ can be partitioned into $I_{1}, I_{2}, \cdots, I_{q}$ and $2\leq q\leq 14$. If $2\leq q \leq 3$, then by the minimality of $q$, we may assume $q=2$. From the above argument, we suppose $|I_{i}|\geq 2$, $i=1, 2$. If the edges from $I_{1}$ to $I_{2}$ are color $c_{i}$ where $1\leq i \leq s$, then there is a monochromatic $W_{5}$, a contradiction. Suppose that the edges from $I_{1}$ to $I_{2}$ are color $c'$ such that $c'\neq c_{i}$ where $1\leq i \leq s$. For each $I_i$ with $i=1,2$, the subgraph induced by the edges in $I_i$ with color $c'$ is a matching. Then
\beqs
|G| & = & |\{v_{1}, v_{2}, \cdots, v_{s}\}|+|I_{1}|+|I_{2}| \\
~ & = & |\{v_{1}, v_{2}, \cdots, v_{s}\}\cup I_{1}|+|I_{2}| \\
~ & \leq & 2[gr_{k'-1}(K_{3}:W_{5})-1]\\
~ & < & n,
\eeqs
a contradiction.

Suppose $4\leq q \leq 14$. Let $I_{1},I_{2},\cdots,I_{r}$ be the parts such that $|I_{i}|\geq 2$ for each $i$ with $1\leq i \leq r$, and $|I_{j}|=1$ for each $j$ with $r+1\leq j \leq q$.
\begin{fact}
$r\leq 5$.
\end{fact}
If $r=5$, then $q=5$ and the reduced graph on the parts $I_1,I_2,I_3,I_4,I_5$ must be the unique $2$-coloring of $K_{5}$ with no monochromatic triangle, say with $I_1I_2I_3I_4I_5I_1$ and $I_1I_3I_5I_2I_4I_1$ making two monochromatic cycles in red and blue respectively. Note that red and blue is not same as $c_i \ (1\leq i\leq s)$. For each $I_i \ (1\leq i\leq 5)$, the subgraph induced by red or blue edges is a matching. For each $I_i \ (1\leq i\leq 5)$, $|\{v_{1}, v_{2}, \cdots, v_{s}\}\cup I_i|\leq gr_{k'-2}(K_{3}:W_{5})-1$, and hence
$$
|G|=|\{v_{1}, v_{2}, \cdots, v_{s}\}|+\sum_{i=1}^5|I_{i}|\leq 5[gr_{k'-2}(K_{3}:W_{5})-1]<n,
$$
a contradiction.

Suppose $r=4$. If $q=4$, then $|\{v_{1}, v_{2}, \cdots, v_{s}\}\cup I_i|\leq gr_{k'-2}(K_{3}:W_{5})-1$ for each $I_i \ (1\leq i\leq 4)$, and hence $|G|=|\{v_{1}, v_{2}, \cdots, v_{s}\}|+\sum_{i=1}^4|I_{i}|\leq 4[gr_{k'-2}(K_{3}:W_{5})-1]<n$, a contradiction. If $q=5$, then $|I_5|=1$ and $|\{v_{1}, v_{2}, \cdots, v_{s}\}\cup I_i|\leq gr_{k'-2}(K_{3}:W_{5})-1$ for each $I_i \ (1\leq i\leq 4)$, and hence $|G|=|\{v_{1}, v_{2}, \cdots, v_{s}\}|+1+\sum_{i=1}^4|I_{i}|\leq 1+4[gr_{k'-2}(K_{3}:W_{5})-1]<n$, a contradiction.

Suppose $r=3$. The triangle in the reduced graph cannot be monochromatic so without loss of generality, suppose all edges from $I_{1}$ to $I_{2} \cup I_{3}$ are red, and $I_2I_3$ is blue. For each $I_i \ (4\leq i\leq q)$, the edges from $I_1$ to $I_i$ is blue.

\begin{claim}\label{Claim3}
$q\leq 7$.
\end{claim}
\begin{proof}
Assume, to the contrary, that $q\geq 8$. Then there are at least five isolated vertices outside $I_{1}\cup I_{2} \cup I_{3}$. Then the subgraph induced by the blue edges in $I_{4}\cup I_{5}\cup \ldots \cup I_{q}$ is a matching, and hence $I_{4}\cup I_{5}\cup \ldots \cup I_{q}$ contains a red $W_5$, a contradiction.
\end{proof}

From Claim \ref{Claim3}, $q\leq 7$. If $5\leq q\leq 7$, then
$|G|=|\{v_{1}, v_{2}, \cdots, v_{s}\}|+4+\sum_{i=1}^3|I_{i}|\leq
4+3[gr_{k'-2}(K_{3}:W_{5})-1]<n$, a contradiction. If $q=4$, then
$|G|=|\{v_{1}, v_{2}, \cdots, v_{s}\}\cup I_1|+1+|I_{2}|+|I_{3}|\leq
1+[gr_{k'-1}(K_{3}:W_{5})-1]+2[gr_{k'-2}(K_{3}:W_{5})-1]<n$, a
contradiction.

Suppose $r=2$. Suppose all edges from $I_1$ to $I_2$ are red. Let $A$ be the set of parts with red edges to $I_1$ and blue edges to $I_2$, and $B$ be the set of parts with blue edges to $I_{1}\cup I_{2}$, and $C$ be the set of parts with blue edges to $I_1$ and red edges to $I_2$.
\begin{claim}\label{Claim4}
$|A|\leq 2$ and $|C|\leq 2$.
\end{claim}
\begin{proof}
Assume, to the contrary, that $|A|\geq 3$. Note that all parts in $A$ are small parts and they are isolated vertices, and hence the edges in $A$ are red or blue. Since $|A|\geq 3$, it follows that there is a vertex of red degree $2$ or a vertex of blue degree $2$, that is, there is a red $P_3$ or blue $P_3$ in $A$. If there is a red $P_3$, then we have a red $W_5$ from $P_3$ and the edges from $A$ to $I_1$, a contradiction. If there is a blue $P_3$, then we have a blue $W_5$ from $P_3$ and the edges from $A$ to $I_2$, also a contradiction.
\end{proof}

From Claim~\ref{Claim4}, we have $|A|\leq 2$ and $|C|\leq 2$. To
avoid a red $W_5$, there is at most a red matching in $I_1$ or
$I_2$, and hence $|I_1|\leq gr_{k'-1}(K_{3}:W_{5})-1$. Since
$|C|\leq 2$, it follows that $C$ contains at most one red edge, and
hence $C\cup I_1$ contains at most a red matching. Clearly,
$|\{v_{1}, v_{2}, \cdots, v_{s}\}\cup C\cup I_1|\leq
gr_{k'-1}(K_{3}:W_{5})-1$.

Suppose $|A\cup B|\leq 2$. Note that there is at most a red matching in $I_2$. If $A\cup B$ contains at most a red matching, then $A\cup B\cup T_2$ contains at most a red matching since the edges from $A\cup B$ to $I_2$ are all blue. If $A\cup B$ contains a vertex of red degree $2$, then the we change all red edges to green (a color different to $c_1,c_2,\ldots,c_s$ and red and blue) and $A\cup B\cup I_2$ contains at most a red matching, and hence
$$
|G|=|\{v_{1}, v_{2}, \cdots, v_{s}\}\cup C\cup H_1|+|A\cup B\cup I_2|\leq 2[gr_{k'-1}(K_{3}:W_{5})-1],
$$
a contradiction. If $|A\cup B|\geq 5$, then $A\cup B$ contains a red $W_5$, a contradiction. If $|A\cup B|=3,4$, then $|I_2|\leq gr_{k'-2}(K_{3}:W_{5})-1$ and hence
\beqs
|G| & = & |\{v_{1}, v_{2}, \cdots, v_{s}\}\cup C\cup I_1|+|A\cup B\cup I_2|\\
~ & \leq & [gr_{k'-1}(K_{3}:W_{5})-1]+4+[gr_{k'-2}(K_{3}:W_{5})-1]\\
~ & < & n,
\eeqs
a contradiction.

Suppose $r=1$. Let $A$ be the set of parts with red edges to $I_1$, and $B$ be the set of parts with blue edges to $I_{1}$.
\begin{claim}\label{Claim5}
$|A|\leq 4$ and $|B|\leq 4$.
\end{claim}
\begin{proof}
Assume, to the contrary, that $|A|\geq 5$. Note that all parts in $A$ are small parts and they are isolated vertices, and hence the edges in $A$ are red or blue. If there is a red $P_3$ in $A$, then there is a red $W_5$ by this $P_3$ and edges from $P_3$ to $I_1$, a contradiction. So $A$ contains at most a red matching and hence there is a blue $W_5$ in $A$, a contradiction.
\end{proof}

From Claim~\ref{Claim5}, we have $|A|\leq 4$ and $|B|\leq 4$. Since
$q\geq 4$, it follows that $|A|\geq 2$ or $|B|\geq 2$. If $|A|\geq
2$ and $|B|\geq 2$, then $|\{v_{1}, v_{2}, \cdots, v_{s}\}\cup
I_1|\leq gr_{k'-2}(K_{3}:W_{5})-1$, and hence
$$
|G|=|\{v_{1}, v_{2}, \cdots, v_{s}\}\cup I_1|+|A\cup B|\leq 8+[gr_{k'-2}(K_{3}:W_{5})-1]<n,
$$
a contradiction. We assume $|A|\geq 2$ and $|B|=1$. Then $|\{v_{1}, v_{2}, \cdots, v_{s}\}\cup I_1|\leq gr_{k'-1}(K_{3}:W_{5})-1$, and hence
$$
|G|=|\{v_{1}, v_{2}, \cdots, v_{s}\}\cup I_1|+|A|+|B|\leq [gr_{k'-1}(K_{3}:W_{5})-1]+4+1<n,
$$
a contradiction.

Suppose $r=0$. If $k'\geq 3$, then $q\leq 14$. Then $|G|=|\{v_{1}, v_{2}, \cdots, v_{s}\}|+q\leq s+q<n$, a contradiction.

Suppose that $v_1,v_2,\ldots,v_s$ does not exist. We assume that $t_1=2$, $|H_{1}^1|\geq 2$ and $|H_{2}^1|\geq 2$. Then $|H_{i}^1|\leq gr_{k'-1}(K_{3}:W_{5})-1$ for $i=1,2$, and hence
$$
|G|=|H_{1}^1|+|H_{2}^1|\leq 2[gr_{k'-1}(K_{3}:W_{5})-1]<n,
$$
a contradiction.

Suppose $4\leq t_1\leq 14$. Note that $H_{1}^1,H_{2}^1,\cdots,H_{r}^1$ be the parts such that $|H_{i}^1|\geq 2$ for each $i \ (1\leq i \leq r)$, and $|H_{j}^1|=1$ for each $j \ (r+1\leq j \leq t_1)$.
\begin{fact}\label{Clm:6}
$r\leq 5$.
\end{fact}

If $r=5$, then $t_1=5$ and the reduced graph on the parts $H_1^1,H_2^1,H_3^1,H_4^1,H_5^1$ must be the unique $2$-coloring of $K_{5}$ with no monochromatic triangle, say with $H_1^1H_2^1H_3^1H_4^1H_5^1H_1^1$ and $H_1^1H_3^1H_5^1H_2^1H_4^1H_1^1$ making two monochromatic cycles in red and blue respectively. For each $H_i^1 \ (1\leq i\leq 5)$, $|H_i^1|\leq gr_{k'-2}(K_{3}:W_{5})-1$, and hence $|G|=\sum_{i=1}^5|H_{i}^1|\leq 5[gr_{k'-2}(K_{3}:W_{5})-1]<n$, a contradiction.

Suppose $r=4$. If $t_1=4$, then $|H_i^1|\leq gr_{k'-2}(K_{3}:W_{5})-1$ for each $H_i^1 \ (1\leq i\leq 4)$, and hence $|G|=\sum_{i=1}^4|H_{i}^1|\leq 4[gr_{k'-2}(K_{3}:W_{5})-1]<n$, a contradiction. If $t_1=5$, then $|H_5^1|=1$ and $|H_i^1|\leq gr_{k'-2}(K_{3}:W_{5})-1$ for each $H_i^1 \ (1\leq i\leq 4)$, and hence $|G|=1+\sum_{i=1}^4|H_{i}^1|\leq 1+4[gr_{k'-2}(K_{3}:W_{5})-1]<n$, a contradiction.

Suppose $r=3$. The triangle in the reduced graph cannot be
monochromatic so without loss of generality, suppose all edges from
$H_{1}^1$ to $H_{2}^1 \cup H_{3}^1$ are red, and $H_2^1H_3^1$ is
blue. For each $H_i^1 \ (4\leq i\leq q)$, the edges from $H_1^1$ to
$H_i^1$ is blue.

\begin{claim}\label{Claim4b}
$t_1\leq 7$.
\end{claim}
\begin{proof}
Assume, to the contrary, that $t_1\geq 8$. Then there are at least five isolated vertices outside $H_{1}^1\cup H_{2}^1 \cup H_{3}^1$. Then the subgraph induced by the blue edges in $H_{4}^1\cup H_{5}^1\cup \ldots \cup H_{t_1}^1$ is a matching, and hence $H_{4}\cup H_{5}^1\cup \ldots \cup H_{t_1}^1$ contains a red $W_5$, a contradiction.
\end{proof}

From Claim~\ref{Claim4b}, $t_1\leq 7$. Then $|G|=4+\sum_{i=1}^3|H_{i}^1|\leq 4+3[gr_{k'-2}(K_{3}:W_{5})-1]<n$, a contradiction.

Suppose $r=2$. Suppose all edges from $H_1^1$ to $H_2^1$ are red. Let $A$ be the set of parts with red edges to $H_1^1$ and blue edges to $H_2^1$, and $B$ be the set of parts with blue edges to $H_{1}^1\cup H_{2}^1$, and $C$ be the set of parts with blue edges to $H_1^1$ and red edges to $H_2^1$.
\begin{fact}
$|A|\leq 2$ and $|C|\leq 2$.
\end{fact}
Clearly, $|C\cup H_1^1|\leq gr_{k'-1}(K_{3}:W_{5})-1$. If $|A\cup B|\leq 2$, then
$$
|G|=|C\cup H_1^1|+|A\cup B\cup H_2^1|\leq 2[gr_{k'-1}(K_{3}:W_{5})-1],
$$
a contradiction. If $|A\cup B|\geq 5$, then $A\cup B$ contains a red $W_5$, a contradiction. If $|A\cup B|=3,4$, then $|H_2^1|\leq gr_{k'-2}(K_{3}:W_{5})-1$ and hence
\beqs
|G| & = & |C\cup H_1|+|A\cup B\cup H_2^1|\\
~ & \leq & [gr_{k'-1}(K_{3}:W_{5})-1]+4+[gr_{k'-2}(K_{3}:W_{5})-1]\\
~ & < & n,
\eeqs
a contradiction.

Suppose $r=1$. Let $A$ be the set of parts with red edges to $H_1^1$, and $B$ be the set of parts with blue edges to $H_{1}^1$. Then $|A|\leq 4$ and $|B|\leq 4$. Then $|A|\geq 2$ or $|B|\geq 2$. If $|A|\geq 2$ and $|B|\geq 2$, then $|H_1^1|\leq gr_{k'-2}(K_{3}:W_{5})-1$, and hence
$$
|G|=|H_1^1|+|A\cup B|\leq 8+[gr_{k'-2}(K_{3}:W_{5})-1]<n,
$$
a contradiction. We assume $|A|\geq 2$ and $|B|=1$. Then $|H_1^1|\leq gr_{k'-1}(K_{3}:W_{5})-1$, and hence
$$
|G|=|H_1^1|+|A|+|B|\leq [gr_{k'-1}(K_{3}:W_{5})-1]+4+1<n,
$$
a contradiction.

Suppose $r=0$. If $k'\geq 3$, then $t_1\leq 14$. Then $|G|=t_1<n$, a contradiction.
\end{proof}

\section{Bounds on the Gallai-Ramsey number For general $n$} \label{Sec:Gen}

For the lower bound, we have the following easy result. We state this result without proof since it follows from the same argument as the proof of Lemma~\ref{Lemma:W5Lower}, in which the value of $R(W_{5}, W_{5})$ is replaced by the lower bounds on $R(W_{n}, W_{n})$ from Lemma~\ref{Lemma:RamseyWnBounds}.

\begin{theorem}
For $k \geq 2$ and $n \geq 6$, we have
$$
gr_{k}(K_{3} : W_{n}) \geq \begin{cases}
(3n - 4) 5^{\frac{k - 2}{2}} + 1 & \text{ if $n$ is even and $k$ is even;}\\
(6n - 8) 5^{\frac{k - 3}{2}} + 1 & \text{ if $n$ is even and $k$ is odd;}\\
(2n - 3) 5^{\frac{k - 2}{2}} + 1 & \text{ if $n$ is odd and $k$ is even;}\\
(4n - 6) 5^{\frac{k - 3}{2}} + 1 & \text{ if $n$ is odd and $k$ is odd.}
\end{cases}
$$
\end{theorem}

We also obtain a general upper bound.

\begin{theorem}\labelz{Thm:GenUpBnd}
For $k \geq 3$ and $n \geq 6$, we have
$$
gr_{k}(K_{3} : W_{n}) \leq (n - 4)^{2} \cdot 30^{k} + k(n - 1).
$$
\end{theorem}

Given nonnegative integers $k, n, r, s, t$ with $k \geq 1$, $n \geq 6$ and $r + s + t = k$, define the number
$$
gr_{k}(K_{3} : rW_{n}, sC_{n - 1}, tP_{n - 2})
$$
to be the minimum integer $N$ such that every $k$-coloring of $K_{N}$ contains one of: a rainbow triangle, a monochromatic copy of $W_{n}$ in one of the first $r$ colors, a monochromatic copy of $C_{n - 1}$ in one of the next $s$ colors, or a monochromatic copy of $P_{n - 2}$ in one of the remaining $t$ colors. In order to prove Theorem~\ref{Thm:GenUpBnd}, we prove the following bound.

\begin{theorem}\labelz{Thm:RefinedGenUp}
Given nonnegative integers $k, n, r, s, t$ with $k \geq 1$, $n \geq 6$ and $r + s + t = k$, we have
$$
gr_{k}(K_{3} : rW_{n}, sC_{n - 1}, tP_{n - 2}) \leq (n - 4)^{2} \cdot 30^{r} \cdot 10^{s} \cdot 2^{t} + k(n - 1).
$$
\end{theorem}

\begin{proof}
Let $G$ be a $k$-coloring of a complete graph of order
$$
N = N(n, r, s, t) = (n - 4)^{2} \cdot 30^{r} \cdot 10^{s} \cdot 2^{t} + k(n + 1)
$$
and suppose that $G$ contains no rainbow triangle, no monochromatic copy of $W_{n}$ in one of the first $r$ colors, no monochromatic copy of $C_{n - 1}$ in one of the remaining $k - r$ colors.

For a colored complete graph $G'$, let $T_{G'}$ be a maximal set of vertices in $G'$ each of which has all one color on its edges to $G' \setminus T_{G'}$ and let $T_{G'}^{i} \subseteq T_{G'}$ be the subset of vertices with all edges of color $i$ to $G' \setminus T_{G'}$ with the additional restriction that $|T_{G'}^{i}| \leq n + 1$ for all $i$. This set of vertices will be called the \emph{garbage set} and vertices will be added to the garbage set only in the process of looking for a monochromatic cycle. In order to avoid creating a monochromatic copy of $C_{n - 1}$ in $G'$, if there is a set $T_{G'}^{i}$ with $|T_{G'}^{i}| \geq \frac{n'}{2}$, then there are no edges of color $i$ within $G' \setminus T_{G'}$ and furthermore, if $n - 1$ is even, there is already a monochromatic copy of $C_{n - 1}$ in color $i$ in $G'$. Since the garbage set always contains at most $n + 1$ vertices corresponding to each color, there will never be more than $k(n + 1)$ vertices in the garbage set, hence the last term in the definition of $N$ above.

Consider a Gallai partition of $G$ with the smallest number of parts $q$ and suppose red is one of the colors that appears on edges in between the parts and blue is the other (if there is a second color). Let $H_{1}, H_{2}, \dots, H_{q}$ be the parts of this partition in decreasing order by their number of vertices. Note that $q \leq R(W_{n}, W_{n}) - 1$. The proof is broken into two main cases based on the parity of $n$.

\setcounter{case}{0}
\begin{case}
$n$ is odd.
\end{case}

Call a part $H_{i}$ of the Gallai partition ``large'' if it has order at least $\frac{n - 1}{2}$. We consider subcases based on the desired red and blue structures.

\begin{subcase} \label{SubCase:TwoWheels}
Both red and blue appear in the first $r$ colors.
\end{subcase}

In this case, we are looking for a red or blue copy of $W_{n}$ in $G$.

First suppose $q \leq 3$ so by the minimality of $q$, we have $q = 2$. Then only red appears in this partition on all edges between $H_{1}$ and $H_{2}$. Then each part $H_{i}$ contains no red copy of $C_{n - 1}$ so
$$
|G| \leq 2[N(n, r - 1, s + 1, t) - 1] < N(n, r, s, t),
$$
a contradiction.

Next suppose $q \geq 4$ so by the minimality of $q$, every part has edges to other parts of the partition in both red and blue. There are at most $5$ large parts of the partition since there can be no monochromatic triangle in the reduced graph among these large parts. In fact, if there are $5$ such parts, then $q = 5$ since any $6$ parts containing $5$ large parts would contain a monochromatic triangle using at least two large parts, yielding a monochromatic copy of $W_{n}$. Thus, there are either at most $5$ parts total or at most $4$ large parts. Since $q \leq R(W_{n}, W_{n})$, we get
\beqs
|G| & = & \sum_{i = 1}^{q} |H_{i}|\\
~ & \leq & \max \begin{cases} 5[ N(n, r - 2, s + 2, t) - 1]\\ 4[N(n, r - 2, s + 2, t) - 1] + [R(W_{n}, W_{n}) - 5] \frac{n - 2}{2} \end{cases} \\
~ & < & N(n, r, s, t),
\eeqs
for a contradiction.

\begin{subcase} \label{SubCase:TwoEvenCycles}
Both red and blue appear in the latter $s + t = k - r$ colors.
\end{subcase}

In this case, we are looking for a red or blue copy of the even cycle $C_{n - 1}$ or path $P_{n - 2}$ in $G$. Since a monochromatic copy of $C_{n - 1}$ contains a monochromatic copy of $P_{n - 2}$, it suffices to find only a monochromatic copy of $C_{n - 1}$.

First suppose $q \leq 3$ so by the minimality of $q$, we have $q = 2$. Then only red appears in this partition on all edges between $H_{1}$ and $H_{2}$. If $|H_{2}| < \frac{n - 1}{2}$, then $H_{2}$ can be added to the garbage set $T_{G}$, contradicting the maximality of $T_{G}$. If $|H_{2}| \geq \frac{n - 1}{2}$, then there is a red copy of $C_{n - 1}$ on the edges between $H_{1}$ and $H_{2}$, for a contradiction.

Next suppose $q \geq 4$ and by minimality of $q$, every part has edges to other parts of the partition in both red and blue. There is at most one large part of the partition to avoid creating a monochromatic copy of $C_{n - 1}$. If $|H_{1}| \geq \frac{n - 1}{2}$, then there are at most $\frac{n - 3}{2}$ vertices in $G \setminus H_{1}$ with red (or similarly blue) edges to $H_{1}$ for a total of at most $n - 3$ vertices in $G \setminus H_{1}$. All of these vertices can be added to $T_{G}$, contradicting the maximality of $T_{G}$. This means that all parts must have order at most $\frac{n - 3}{2}$. With at most $R(C_{n - 1}, C_{n - 1}) = n - 2 + \frac{n - 1}{2}$ parts, this means that
$$
|G| \leq \frac{n - 3}{2}\left[n - 2 + \frac{n - 1}{2}\right],
$$
a contradiction.

\begin{subcase}
One of red or blue (say red) appears in the first $r$ colors while the other appears among the latter $s + t = k - r$ colors.
\end{subcase}

In this case, we are looking for a red copy of $W_{n}$ or a blue copy of $C_{n - 1}$ or $P_{n - 2}$ in $G$. Since a blue copy of $C_{n - 1}$ contains a blue copy of $P_{n - 2}$, it suffices to find only a blue copy of $C_{n - 1}$.

First suppose $q \leq 3$ so by the minimality of $q$, we have $q = 2$. Then only one color appears on edges between the two parts $H_{1}$ and $H_{2}$ and we may apply one of the previous two subcases.

Next suppose $q \geq 4$ and by minimality of $q$, every part has edges to other parts of the partition in both red and blue. There are at most $2$ large parts of the partition since there can be no red triangle in the reduced graph among the large parts and no blue edge in the reduced graph among the large parts.

If $|H_{1}| \geq \frac{n - 1}{2}$, then there are at most $\frac{n - 3}{2}$ vertices with blue edges to $H_{1}$, call that set $B$ and the set of vertices remaining in $G \setminus (H_{1} \cup B)$ is called $A$. Then $H_{1}$ and $A$ each contain no red copy of $C_{n - 1}$ so
$$
|G| \leq \frac{n - 3}{2} + 2[N(n, r - 1, s + 1, t) - 1] < N(n, r, s, t),
$$
a contradiction.

\begin{case}
$n$ is even.
\end{case}

In this case, we call a part $H_{i}$ of the Gallai partition ``large'' if it has order at least $\frac{n - 2}{2}$.

\begin{subcase}
Both red and blue appear in the first $r$ colors.
\end{subcase}

This subcase follows exactly the same argument as Subcase~\ref{SubCase:TwoWheels}.

\begin{subcase} 
Both red and blue appear in the middle $s$ colors.
\end{subcase}

In order to avoid a red or blue copy of $C_{n - 1}$, we must have $q \leq [2(n - 1) - 1] - 1 = 2n - 6$.

If $q \leq 3$, then by minimality of $q$, we may assume $q = 2$, say with red edges appearing in between the two parts. Then if either part is large, it contains no red copy of $P_{n - 2}$, so we have
\beqs
|G| & = & |H_{1}| + |H_{2}|\\
~ & \leq & 2 [N(n, r, s - 1, t + 1) - 1]\\
~ & < & N(n, r, s, t),
\eeqs
a contradiction.

Thus, suppose $q \geq 4$ and by minimality of $q$, each part has edges to other parts in both red and blue. Since a monochromatic triangle in the reduced graph restricted to large parts would contain a monochromatic copy of $C_{n - 1}$, there can be at most $5$ large parts. Each of these large parts contains no red or blue path $P_{n - 2}$ so we have
\beqs
|G| & = & \sum_{i = 1}^{q} |H_{i}|\\
~ & \leq & 5 [N(n, r, s - 1, t + 1) - 1] + [(2n - 6) - 5]\left[ \frac{n - 4}{2} \right]\\
~ & < & N(n, r, s, t),
\eeqs
a contradiction.

\begin{subcase} 
Both red and blue appear in the last $t$ colors.
\end{subcase}

This subcase follows exactly the same argument as Subcase~\ref{SubCase:TwoEvenCycles}.

Note that for the remaining subcases, we may assume $q \geq 4$ since otherwise the proof reduces to one of the first three subcases. By minimality of $q$, each part has edges to some other parts in red and some others in blue.

\begin{subcase} 
Red appears in the first $r$ colors and blue appears in the next $s$ colors.
\end{subcase}

In this case, we have $4 \leq q \leq R(W_{n}, C_{n - 1}) - 1 \leq 3n - 2$ (see \cite{CCMN09} for example). In order to avoid a red copy of $W_{n}$ or a blue copy of $C_{n - 1}$, there can be at most $R(K_{4}, K_{3}) - 1 = 8$ large parts. Each of these large parts contains no red copy of $C_{n - 1}$ and no blue copy of $P_{n - 2}$. This means that
\beqs
|G| & = & \sum_{i = 1}^{q} |H_{i}|\\
~ & \leq & 8[N(n, r - 1, s, t + 1) - 1] + [(3n - 2) - 8]\left[ \frac{n - 4}{2} \right]\\
~ & < & N(n, r, s, t),
\eeqs
a contradiction.

\begin{subcase} 
Red appears in the first $r$ colors and blue appears in the last $t$ colors.
\end{subcase}

In this case, we have $4 \leq q \leq R(W_{n}, P_{n - 2}) - 1 \leq 3n - 2$ (using the same results as cited above loosely).  In order to avoid a red copy of $W_{n}$ or a blue copy of $P_{n - 2}$, there can be at most $3$ large parts and in between these large parts must only be red edges. Each of these large parts contains no red copy of $C_{n - 1}$ and no blue copy of $P_{n - 2}$. This means that
\beqs
|G| & = & \sum_{i = 1}^{q} |H_{i}|\\
~ & \leq & 3[N(n, r - 1, s + 1, t) - 1] + [(3n - 2) - 3]\left[ \frac{n - 4}{2} \right]\\
~ & < & N(n, r, s, t),
\eeqs
a contradiction.

\begin{subcase} 
Red appears in the middle $s$ colors and blue appears in the last $t$ colors.
\end{subcase}

In this case, we have $4 \leq q \leq R(C_{n - 1}, P_{n - 2}) - 1 \leq \frac{3(n - 2)}{2} - 1$. In order to avoid a red copy of $C_{n - 1}$ or a blue copy of $P_{n - 2}$, there can be at most $2$ large parts and in between these large parts must only be red edges. Each of these large parts contains no red or blue copy of $P_{n - 2}$. This means that
\beqs
|G| & = & \sum_{i = 1}^{q} |H_{i}|\\
~ & \leq & 2[N(n, r - 1, s + 1, t) - 1] + \left[\frac{3(n - 2)}{2} - 1 - 2\right]\left[ \frac{n - 4}{2} \right]\\
~ & < & N(n, r, s, t),
\eeqs
a contradiction.
\end{proof}


\end{document}